\documentclass[a4paper,12pt]{amsart}

\usepackage{amsthm,amsmath,xspace,url,gensymb}
\newcommand{\Dfn}[1]{{\sf #1}}                            % for definitions

\usepackage[breaklinks=true,hyperref]{hyperref}
\usepackage[line,matrix,frame,curve,poly]{xy}

\newtheorem{thm}{Theorem}[section]
\newtheorem{prop}[thm]{Proposition}
\newtheorem{cor}[thm]{Corollary}
\newtheorem{lem}[thm]{Lemma}

\theoremstyle{definition}

\newtheorem{eg}[thm]{Example}
\newtheorem*{rmk}{Remark}

\def\op{\operatorname{op}}
\def\cl{\operatorname{cl}}
\def\B{\mathcal{B}}

\newcommand{\Set}[1]{\mathcal{#1}}
\newcommand{\abs}[1]{\big\lvert #1\big\rvert}

\def\An{\operatorname{A}_{n-1}}

\begin{document}
\title[Crossings and nestings in set partitions of classical types]{Crossings
  and nestings in set partitions \\ of classical types} 
\author{Martin Rubey}
\address{Institut f\"ur Algebra, Zahlentheorie und Diskrete Mathematik,Leibniz
  Universit\"at Hannover, Welfengarten 1, D-30167 Hannover, Germany}
\email{martin.rubey@math.uni-hannover.de}
\urladdr{http://www.iazd.uni-hannover.de/~rubey/}
        
\author{Christian Stump} 
\address{Fakult\"at f\"ur Mathematik, Universit\"at
  Wien, Nordbergstra{\ss}e 15, A-1090 Vienna, Austria}
\email{christian.stump@univie.ac.at}
\urladdr{http://homepage.univie.ac.at/christian.stump/}

\date{\today}

\begin{abstract}
  In this article, we investigate bijections on various classes of set
  partitions of classical types that preserve openers and closers.  On the one
  hand we present bijections that interchange crossings and nestings.  For
  types $B$ and $C$, they generalize a construction by Kasraoui and Zeng for
  type $A$, whereas for type $D$, we were only able to construct a bijection
  between non-crossing and non-nesting set partitions.  On the other hand we
  generalize a bijection to type $B$ and $C$ that interchanges the cardinality
  of the maximal crossing with the cardinality of the maximal nesting, as given
  by Chen, Deng, Du, Stanley and Yan for type $A$.  Using a variant of this
  bijection, we also settle a conjecture by Soll and Welker concerning
  generalized type $B$ triangulations and symmetric fans of Dyck paths.
\end{abstract}

\maketitle
        
\setcounter{tocdepth}{2}
\tableofcontents

%%%%%%%%%%%%%%%%%%%%%%%%%%%%%%%%%%%%%%%%%%%%%%%%%%%%%%%%%%%%%%%%%%%%%%%%%%%%%%%
\section*{Introduction}
%%%%%%%%%%%%%%%%%%%%%%%%%%%%%%%%%%%%%%%%%%%%%%%%%%%%%%%%%%%%%%%%%%%%%%%%%%%%%%%

The lattice of non-crossing set partitions was first considered by Germain
Kreweras in \cite{Kreweras1972}.  It was later reinterpreted by Paul Edelman,
Rodica Simion and Daniel Ullman, as a well-behaved sub-lattice of the
intersection lattice for the hyperplane arrangement of type $A$, see
e.g. \cite{Edelman1980, EdelmanSimion1994, SimionUllman1991}.  Natural
combinatorial interpretations of non-crossing partitions for the classical
reflection groups were then given by Christos A.~Athanasiadis and Vic Reiner 
in~\cite{Reiner1997,AthanasiadisReiner2004}.
                
On the other hand, non-nesting partitions were simultaneously introduced for
all crystallographic reflection groups by Alex Postnikov as anti-chains in the
associated root poset, see \cite[Remark 2]{Reiner1997}.
                
Within the last years, several bijections between non-crossing and non-nesting
partitions have been constructed.  In particular, type (i.e., block-size)
preserving bijections were given by Christos
A.~Athanasiadis~\cite{Athanasiadis1998} for type $A$ and by Alex Fink and
Benjamin I.~Giraldo~\cite{FinkGiraldo2009} for the other classical reflection
groups.  One of the authors of the present article~\cite{Stump2008} constructed
another bijection for types $A$ and $B$ which transports other natural
statistics.  Recently, Ricardo Mamede~\cite{Mamede2009} constructed a bijection
for types $A$ and $B$ which turns out to be subsumed by the bijections we
will present here.

The material on non-crossing partitions on the one hand and on non-nesting
partitions on the other hand suggests that they are not only counted by the
same numbers, namely the Catalan numbers, but
are more deeply connected.  These connections were explored by Drew Armstrong
in \cite[Chapter 5.1.3]{Armstrong2007}.  In this paper we would like to exhibit
some further connections.
                
In the case of set partitions of type $A$, also the \emph{number} of crossings
and nestings was considered: Anisse Kasraoui and Jiang Zeng constructed a
bijection which interchanges crossings and nestings in \cite{KasraouiZeng2006}.
Finally, in a rather different direction, William~Y.C. Chen, Eva~Y.P. Deng,
Rosena~R.X. Du, Richard~P. Stanley~\cite{ChenDengDuStanleyYan2006} have shown
that the number of set partitions where the \emph{maximal crossing} has
cardinality $k$ and the \emph{maximal nesting} has cardinality $l$ is the same
as the number of set partitions where the maximal crossing has cardinality $l$
and the maximal nesting has cardinality $k$.
                
In this paper, we present bijections on various classes of set partitions of
classical types that preserve openers and closers.  In particular, the
bijection by Anisse Kasraoui and Jiang Zeng as well as the bijection by
William~Y.C. Chen, Eva~Y.P. Deng, Rosena~R.X. Du, Richard~P. Stanley enjoy this
property.  We give generalizations of these bijections for the other classical
reflection groups, whenever possible.  Furthermore we show that the bijection
is in fact (mostly) unique for the class of non-crossing and non-nesting set
partitions.  Finally, a slight variation of one of our bijections settles a
conjecture by Daniel Soll and Volkmar Welker~\cite{SollWelker2006}, concerning
generalized triangulations with $180\degree$ symmetry and symmetric fans of
Dyck paths.

%%%%%%%%%%%%%%%%%%%%%%%%%%%%%%%%%%%%%%%%%%%%%%%%%%%%%%%%%%%%%%%%%%%%%%%%%%%%%%%
\section{Set partitions for classical types}
%%%%%%%%%%%%%%%%%%%%%%%%%%%%%%%%%%%%%%%%%%%%%%%%%%%%%%%%%%%%%%%%%%%%%%%%%%%%%%%
\label{backgroundanddefinitions}

A \Dfn{set partition} of $[n]:= \{1,2,3,\ldots,n\}$ is a collection $\B$ of
pairwise disjoint, non-empty subsets of $[n]$, called \Dfn{blocks}, whose union
is $[n]$.  We visualize $\B$ by placing the numbers $1,2,\dots, n$ in this
order on a line and then joining \emph{consecutive} elements of each block by
an arc, see Figures~\ref{NC9} and \ref{NN9} for examples.

\begin{figure}
  \setlength{\unitlength}{25pt}
  \begin{picture}(8,2.5)
    \put(6,0){\hbox{$7$}}
    \put(8,0){\hbox{$9$}}
    \qbezier(6.2,0.5)(7.15,1.5)(8.1,0.5)
    \put(0,0){\hbox{$1$}}
    \put(1,0){\hbox{$2$}}
    \put(2,0){\hbox{$3$}}
    \put(3,0){\hbox{$4$}}
    \put(4,0){\hbox{$5$}}
    \put(5,0){\hbox{$6$}}
    \put(7,0){\hbox{$8$}}
    \qbezier(0.2,0.5)(3.15,3.5)(6.1,0.5)
    \qbezier(1.2,0.5)(2.65,2.5)(4.1,0.5)
    \qbezier(4.2,0.5)(4.65,1.5)(5.1,0.5)
    \qbezier(2.2,0.5)(2.65,1.5)(3.1,0.5)
  \end{picture}
  \caption{A non-crossing set partition of $[9]$.}
  \label{NC9}
\end{figure}
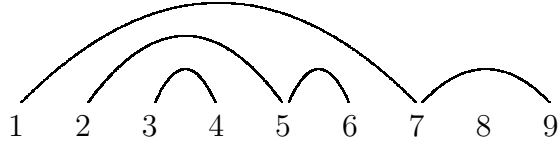
\begin{figure}
  \setlength{\unitlength}{25pt}
  \begin{picture}(8,1.5)
    \put(0,0){\hbox{$1$}}
    \put(1,0){\hbox{$2$}}
    \put(2,0){\hbox{$3$}}
    \put(3,0){\hbox{$4$}}
    \put(4,0){\hbox{$5$}}
    \put(5,0){\hbox{$6$}}
    \put(6,0){\hbox{$7$}}
    \put(7,0){\hbox{$8$}}
    \put(8,0){\hbox{$9$}}
    \qbezier(0.2,0.5)(1.65,2.5)(3.1,0.5)
    \qbezier(1.2,0.5)(2.65,2.5)(4.1,0.5)
    \qbezier(2.2,0.5)(3.65,2.5)(5.1,0.5)
    \qbezier(4.2,0.5)(5.15,1.5)(6.1,0.5)
    \qbezier(6.2,0.5)(7.15,1.5)(8.1,0.5)
  \end{picture}
  \caption{A non-nesting set partition of $[9]$.}
  \label{NN9}
\end{figure}
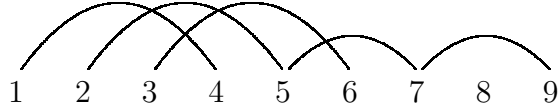	

The \Dfn{openers} $\op(\B)$ are the non-maximal elements of the blocks in $\B$,
whereas the \Dfn{closers} $\cl(\B)$ are its non-minimal elements. For example,
the set partitions in Figures~\ref{NC9} and \ref{NN9} both have $\op(\B) =
\{1,2,3,5,7\}$ and $\cl(\B) = \{4,5,6,7,9\}$.

A pair $(\Set{O},\Set{C})\subseteq [n]\times [n]$ is an \Dfn{opener-closer
  configuration}, if $\abs{\Set{O}}=\abs{\Set{C}}$ and
$$
\abs{\Set{O}\cap [k]}\geq\abs{\Set{C}\cap [k+1]} \quad\text{for}\quad
k\in\{0,1,\dots,n-1\},
$$
or, equivalently, $(\Set{O},\Set{C})=\big(\op(\B),\cl(\B)\big)$ for some set
partition $\B$ of $n$.

We remark that in~\cite{KasraouiZeng2006}, Anisse Kasraoui and Jiang Zeng
distinguish between openers, closers and \emph{transients}, which are, in our
definition, those numbers which are both openers and closers.

It is now well established that set partitions of $[n]$ are in natural
bijection with intersections of the reflecting hyperplanes $x_i-x_j=0$ in
$\mathbb R^n$ of the Coxeter group of type $\An$.  For example, the set
partition in Figure~\ref{NC9} corresponds to the intersection
$$\{x\in\mathbb R^9: x_1=x_7=x_9, x_2=x_5=x_6, x_3=x_4\}.$$

Therefore, set partitions of $[n]$ can be seen as set partitions of type $\An$
and set partitions of other types can be defined by analogy,
see~\cite{Athanasiadis1998, Reiner1997}.  The reflecting hyperplanes for $B_n$
and $C_n$ are
\begin{itemize}
\item[] $x_i=0$ for $1\leq i \leq n$,
\item[] $x_i-x_j=0$ for $1\leq i < j \leq n$, and
\item[] $x_i+x_j=0$ for $1\leq i < j \leq n$.
\end{itemize}
Thus, a set partition of type $B_n$ or $C_n$ is a set partition $\B$ of $[\pm
n] := \{1,2,\ldots,n,-1,-2,\ldots,-n\}$, such that
\begin{eqnarray}
  B \in \B \Leftrightarrow -B \in \B \label{eq:setpartitionsymmetry}
\end{eqnarray}
and such that there exists at most one block $B_0 \in \B$ (called the \Dfn{zero
  block}) for which $B_0 = -B_0$.

The hyperplanes for $D_n$ are those for $B_n$ and $C_n$ other than $x_i=0$ for
$1\leq i \leq n$, whence a set partition $\B$ of type $D_n$ is a set partition
of type $B_n$ (or $C_n$) where the zero block, if present, must not consist of
a single pair $\{i, -i\}$.

%%%%%%%%%%%%%%%%%%%%%%%%%%%%%%%%%%%%%%%%%%%%%%%%%%%%%%%%%%%%%%%%%%%%%%%%%%%%%%%
\section{Crossings and nestings in set partitions of type $A$}
%%%%%%%%%%%%%%%%%%%%%%%%%%%%%%%%%%%%%%%%%%%%%%%%%%%%%%%%%%%%%%%%%%%%%%%%%%%%%%%
\label{sec:cross-nest-A}

One of the goals of this article is to refine the following well known correspondences between
non-crossing and non-nesting set partitions.  For ordinary set partitions, a
\Dfn{crossing} consists of a pair of arcs $(i,j)$ and $(i',j')$ such that
$i<i'<j<j'$:
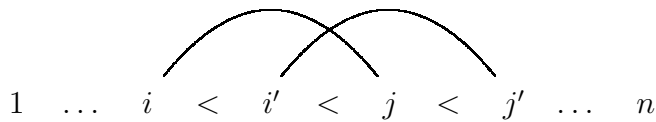
\begin{center}
  \setlength{\unitlength}{20pt}
  \begin{picture}(12,2.2)
    \put(0,0){\hbox{$1$}}
    \put(1,0){\hbox{$\ldots$}}
    \put(2.5,0){\hbox{$i$}}
    \put(3.5,0){\hbox{$<$}}
    \put(4.75,0){\hbox{$i'$}}
    \put(5.75,0){\hbox{$<$}}
    \put(7,0){\hbox{$j$}}
    \put(8,0){\hbox{$<$}}
    \put(9.25,0){\hbox{$j'$}}
    \put(10.25,0){\hbox{$\ldots$}}
    \put(11.75,0){\hbox{$n$}}
    \qbezier(2.9,0.7)(4.9,3.2)(6.9,0.7)
    \qbezier(5.1,0.7)(7.1,3.2)(9.1,0.7)
  \end{picture}
\end{center}
On the other hand, if $i<i'<j'<j$, we have a \Dfn{nesting}, pictorially:
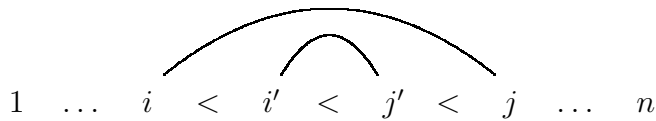
\begin{center}
  \setlength{\unitlength}{20pt}
  \begin{picture}(12,2.2)
    \put(0,0){\hbox{$1$}}
    \put(1,0){\hbox{$\ldots$}}
    \put(2.5,0){\hbox{$i$}}
    \put(3.5,0){\hbox{$<$}}
    \put(4.75,0){\hbox{$i'$}}
    \put(5.75,0){\hbox{$<$}}
    \put(7,0){\hbox{$j'$}}
    \put(8,0){\hbox{$<$}}
    \put(9.25,0){\hbox{$j$}}
    \put(10.25,0){\hbox{$\ldots$}}
    \put(11.75,0){\hbox{$n$}}
    \qbezier(2.9,0.7)(6.0,3.2)(9.1,0.7)
    \qbezier(5.1,0.7)(6.0,2.2)(6.9,0.7)
  \end{picture}
\end{center}
A set partition of $[n]$ is called \Dfn{non-crossing} (resp. \Dfn{non-nesting})
if the number of crossings (resp. the number of nestings) equals $0$.

It has been known for a long time that the numbers of non-crossing and
non-nesting set-partitions of $[n]$ coincide.  More recently, Anisse Kasraoui
and Jiang Zeng have shown in~\cite{KasraouiZeng2006} that much more is true:
\begin{thm}
  There is an explicit bijection on set partitions of $[n]$, preserving the set
  of openers and the set of closers, and interchanging the number of crossings
  and the number of nestings.
\end{thm}
The construction in~\cite{KasraouiZeng2006} also proves the following corollary:
\begin{cor}\label{cor:unique-A}
  For any opener-closer configuration $(\Set{O},\Set{C}) \subseteq [n] \times
  [n]$, there exists a unique non-crossing set partition $\B$ of $[n]$ and a
  unique non-nesting set partition $\B'$ of $[n]$ such that
  $$
  \op(\B)=\op(\B')=\Set{O} \quad\text{and}\quad \cl(\B)=\cl(\B')=\Set{C}.
  $$
\end{cor}

In the following section, we will provide a proof completely analogous to the
one of Anisse Kasraoui and Jiang Zeng, for Type $C$.

Apart from the number of crossings or nestings, another natural statistic to
consider is the cardinality of a \lq maximal crossing\rq\ and of a \lq maximal
nesting\rq: a \Dfn{maximal crossing} of a set partition is a set of largest
cardinality of mutually crossing arcs, whereas a \Dfn{maximal nesting} is a set
of largest cardinality of mutually nesting arcs.  For example, in
Figure~\ref{NC9}, the arcs $\{(1,7),(2,5),(3,4)\}$ are a maximal nesting of
cardinality $3$.  In Figure~\ref{NN9} the arcs $\{(1,4),(2,5),(3,6)\}$ are a
maximal crossing.

The following symmetry property was shown by William~Y.C. Chen, Eva~Y.P. Deng,
Rosena~R.X. Du, Richard~P. Stanley and
Catherine~H. Yan~\cite{ChenDengDuStanleyYan2006}:
\begin{thm}\label{thm:k-crossing-A}
  There is an explicit bijection on set partitions, preserving the set of
  openers and the set of closers, and interchanging the cardinalities of the
  maximal crossing and the maximal nesting.
\end{thm}

Since a \lq maximal crossing\rq\ of a non-crossing partition and a \lq maximal
nesting\rq\ of a non-nesting partition both have cardinality $1$,
Corollary~\ref{cor:unique-A} implies that this bijection coincides with the
bijection by Anisse Kasraoui and Jiang Zeng for non-crossing and non-nesting
partitions.  In particular, we obtain the curious fact that in this case, the
bijection maps non-crossing partitions with $k$ nestings and maximal nesting
having cardinality $l$ to non-nesting partitions with $k$ crossings and maximal
crossing having cardinality $l$.

We have to stress however, that in general it is not possible to interchange
the number of crossings and the cardinality of a maximal crossing with the
number of nestings and the cardinality of a maximal nesting simultaneously.

\begin{eg}
  For $n=8$, there is a set partition with one crossing, six nestings and the
  cardinalities of the maximal crossing and the maximal nesting equal both one,
  namely $\{\{1,7\},\{2,8\},\{3,4,5,6\}\}$.  However, there is no set partition
  with six crossings, one nesting and cardinalities of the maximal crossing and
  the maximal nesting equal to one.  To check, the four set partitions with six
  crossings and one nesting are 
  \begin{align*}
    &\{\{1,4,6\},\{2,5,8\},\{3,7\}\},\\
    &\{\{1,4,7\},\{3,5,8\},\{2,6\}\},\\
    &\{\{1,4,8\},\{2,5,7\},\{3,6\}\},\\
    &\{\{1,5,8\},\{2,4,7\},\{3,6\}\}.
  \end{align*}
\end{eg}

%%%%%%%%%%%%%%%%%%%%%%%%%%%%%%%%%%%%%%%%%%%%%%%%%%%%%%%%%%%%%%%%%%%%%%%%%%%%%%%
\section{Crossings and nestings in set partitions of type $C$}
%%%%%%%%%%%%%%%%%%%%%%%%%%%%%%%%%%%%%%%%%%%%%%%%%%%%%%%%%%%%%%%%%%%%%%%%%%%%%%%
\label{sec:cross-nest-C}

Type independent definitions for \Dfn{non-crossing} and \Dfn{non-nesting} set
partitions have been available for a while now, see for example
\cite{Armstrong2007, Athanasiadis1998, AthanasiadisReiner2004, Reiner1997}.
However, it turns out that the notions of crossing and nesting is more tricky,
and we do not have a type independent definition.  In this section we
generalize the results of the previous section to type $C$.

We want to associate two pictures to each set partition, namely the
\lq crossing\rq\ and the \lq nesting diagram\rq.  To this end, we define two
orderings on the set $[\pm n]$: the \Dfn{nesting order} for type $C$ is
$$1<2<\dots<n<-n<\dots<-2<-1,$$
whereas the \Dfn{crossing order} is
$$1<2<\dots<n<-1<-2<\dots<-n.$$

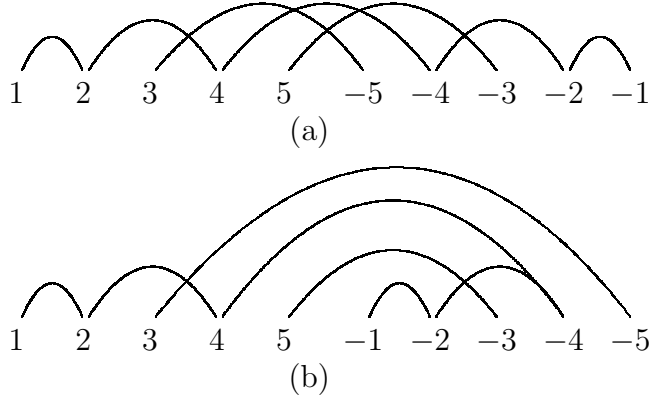
\begin{figure}
  \setlength{\unitlength}{25pt}
  \begin{picture}(9,1.5)
    \put(0,0){\hbox{$1$}}
    \put(1,0){\hbox{$2$}}
    \put(2,0){\hbox{$3$}}
    \put(3,0){\hbox{$4$}}
    \put(4,0){\hbox{$5$}}
    \put(5,0){\hbox{$-5$}}
    \put(6,0){\hbox{$-4$}}
    \put(7,0){\hbox{$-3$}}
    \put(8,0){\hbox{$-2$}}
    \put(9,0){\hbox{$-1$}}
    \qbezier(0.2,0.5)(0.65,1.5)(1.1,0.5)
    \qbezier(1.2,0.5)(2.15,2)(3.1,0.5)
    \qbezier(2.2,0.5)(3.85,2.5)(5.3,0.5)	
    \qbezier(3.2,0.5)(4.75,2.5)(6.3,0.5)
    \qbezier(4.2,0.5)(5.75,2.5)(7.3,0.5)
    \qbezier(6.4,0.5)(7.35,2)(8.3,0.5)
    \qbezier(8.4,0.5)(8.85,1.5)(9.3,0.5)
  \end{picture}\\
  (a) \\
  \setlength{\unitlength}{25pt}
  \begin{picture}(9,3)
    \put(0,0){\hbox{$1$}}
    \put(1,0){\hbox{$2$}}
    \put(2,0){\hbox{$3$}}
    \put(3,0){\hbox{$4$}}
    \put(4,0){\hbox{$5$}}
    \put(5,0){\hbox{$-1$}}
    \put(6,0){\hbox{$-2$}}
    \put(7,0){\hbox{$-3$}}
    \put(8,0){\hbox{$-4$}}
    \put(9,0){\hbox{$-5$}}
    \qbezier(0.2,0.5)(0.65,1.5)(1.1,0.5)
    \qbezier(1.2,0.5)(2.15,2)(3.1,0.5)
    \qbezier(2.2,0.5)(5.85,5)(9.3,0.5)
    \qbezier(3.2,0.5)(5.75,4)(8.3,0.5)
    \qbezier(4.2,0.5)(5.75,2.5)(7.3,0.5)
    \qbezier(6.4,0.5)(7.35,2)  (8.3,0.5)
    \qbezier(5.4,0.5)(5.85,1.5)(6.3,0.5)
  \end{picture}\\
  (b)
  \caption{The nesting (a) and the crossing (b) diagram of a set partition of
    type $C_5$.}
  \label{fig:NNC5}
\end{figure}

The \Dfn{nesting diagram} of a set partition $\B$ of type $C_n$ is obtained by
placing the numbers in $[\pm n]$ in \emph{nesting order} on a line and then
joining consecutive elements of each block of $\B$ by an arc, see
Figure~\ref{fig:NNC5}(a) for an example.

The \Dfn{crossing diagram} of a set partition $\B$ of type $C_n$ is obtained
from the \Dfn{nesting diagram} by reversing the order of the negative numbers.
More precisely, we place the numbers in $[\pm n]$ in \emph{crossing order} on a
line and then join \emph{consecutive elements in the nesting order} of each
block of $\B$ by an arc, see Figure~\ref{fig:NNC5}(b) for an example.  We
stress that the same elements are joined by arcs in both diagrams.  Observe
furthermore that the symmetry property \eqref{eq:setpartitionsymmetry} implies
that if $(i,j)$ is an arc, then its negative $(-j,-i)$ is also an arc.

A \Dfn{crossing} is a pair of arcs that crosses in the crossing diagram, and a
\Dfn{nesting} is a pair of arcs that nests in the nesting diagram.

The \Dfn{openers} $\op(\B)$ are the \emph{positive} non-maximal elements of the
blocks in $\B$ with respect to the \emph{nesting order}, the \Dfn{closers}
$\cl(\B)$ the \emph{positive} non-minimal elements.  Thus, openers and closers
are the start and end points of the arcs in the positive part of the nesting
diagram. For example, the set partition displayed in Figure~\ref{fig:NNC5} has
openers $\{1,2,3,4,5\}$ and closers $\{2,4\}$.

In type $C_n$, $(\Set{O},\Set{C})\subseteq [n]\times [n]$ is an
\Dfn{opener-closer configuration}, if
$$\label{eq:OCconfiguration}
\abs{\Set{O}\cap [k]}\geq\abs{\Set{C}\cap [k+1]}\quad\text{for}\quad
k\in\{0,1,\dots,n-1\}.
$$
Note that we do not require that $\abs{\Set{O}}=\abs{\Set{C}}$.  For
convenience, we call the negatives of the elements in $\Set{O}$ \Dfn{negative
  closers} and the negatives of the elements in $\Set{C}$ \Dfn{negative
  openers}.

\begin{thm}\label{thm:crossings-C}
  There is an explicit bijection on set partitions of type $C$, preserving the
  set of openers and the set of closers, and interchanging the number of
  crossings and the number of nestings.
\end{thm}
In fact, the proof will show that the statement of the theorem remains valid if
we restrict ourselves to arcs that have a positive opener.

Furthermore, we will also see the following analog of
Corollary~\ref{cor:unique-A}:
\begin{cor}
  For any opener-closer configuration $(\Set{O},\Set{C}) \subseteq
  [n] \times [n]$, there exists a unique non-crossing set partition $\B$ and a
  unique non-nesting set partition $\B'$, both of type $C_n$, such that 
  $$
  \op(\B)=\op(\B')=\Set{O} \quad\text{and}\quad \cl(\B)=\cl(\B')=\Set{C}.
  $$
\end{cor}
\begin{proof}
  The bijection proceeds in three steps.  In the first step we consider only
  the given opener-closer configuration, and connect every closer, starting
  with the smallest, with the appropriate opener.  Let us call an opener
  \Dfn{active}, if it has not yet been connected with a closer.

  Let $\B$ be a set partition of type $C_n$.  Every closer in $\B$ corresponds
  to an arc $(i,j)$ with positive $j$ in the given set partition.  It is nested
  by precisely those arcs $(i',j')$ that have opener $1 \leq i'<i$ and closer
  $j<j'\leq n$ or negative closer.  On the other hand, it is crossed by those
  arcs $(i',j')$ that have opener $i'$ with $i<i'<j$ and closer $j<j' \leq n$
  or negative closer.

  To construct the image of $\B$, we want to interchange the number of arcs
  crossing the arc ending in $j$ with the number of arcs nesting it.  Thus, if
  there are $k$ active openers smaller than $j$, and $(i,j)$ is crossed by $c$
  arcs in $\B$, we connect $j$ with the $(c+1)$\textsuperscript{st} active
  opener.  Then, the arc ending in $j$ will be nested by precisely $c$ arcs.
  The first step is completed when all closers have been connected.

  Note that we do not have any choice if we want to construct, say, a
  non-nesting set partition: connecting $j$ with any other active opener but the
  very first will produce a nesting.

  In the second step, we use the symmetry property
  \eqref{eq:setpartitionsymmetry} to connect elements $(i',j')$ with both $i'$
  and $j'$ negative.  More precisely, for every arc $(i,j)$ with (positive)
  closer $j$, we add an arc $(-j, -i)$ to the set partition we are
  constructing.

  Finally, we need to connect the remaining active openers with appropriate
  negative closers.  Observe that two arcs $(i,j)$ and $(i',j')$ where both $i$
  and $i'$ are positive and both $j$ and $j'$ are negative cross if and only if
  they nest.  Suppose that the arcs connecting positive with negative elements
  in $\B$ are $\{(i_1,j_1), (i_2,j_2),\dots,(i_k,j_k)\}$.  Obviously, the set
  $\{i_1, i_2,\dots,i_k\}$ and $\{-j_1, -j_2,\dots,-j_k\}$ are identical, and
  the arcs define a matching $\sigma$, such that $j_m=-i_{\sigma(m)}$.

  Thus, if the remaining active openers are $\{o_1, o_2,\dots,o_k\}$, the image
  of $\B$ shall contain the arcs
  $\{(o_1,-o_{\sigma(1)}),(o_2,-o_{\sigma(2)}),\dots,(o_k,-o_{\sigma(k)})\}$.
  This completes the description of the bijection.

  Again, note that we do not have any choice if we want to construct a
  non-nesting or non-crossing set partition: there is only one non-crossing --
  and therefore only one non-nesting -- matching of the appropriate size that
  satisfies the symmetry property \eqref{eq:setpartitionsymmetry}.
\end{proof}

In Section~\ref{sec:k-crossing-C} we will show the following analog to
Theorem~\ref{thm:k-crossing-A}, where the definition of maximal crossing is as
in type $A$:
\begin{thm}\label{thm:k-crossing-C}
  There is an explicit bijection on set partitions of type $C$, preserving the
  set of openers and the set of closers, and interchanging the cardinalities of
  the maximal crossing and the maximal nesting.
\end{thm}
\begin{rmk}
  It is tempting to consider a different notion of crossing and nesting, as
  suggested by Drew Armstrong in~\cite{Armstrong2007}.  He defined a
  \emph{bump} as an equivalence class of arcs, where the arc $(i,j)$ is
  identified with $(-j, -i)$.  From an algebraic point of view this is a very
  natural idea, since both correspond to the same hyperplane $x_i=x_j$, or, if
  $i$ and $j$ have opposite signs, to $x_i=-x_j$.

  As an example, the partition $\{(1,4,-2),(3,5)\}$ would then be $3$-crossing,
  since with this definition $(1,4)$ crosses $(3,5)$ but also $(2,-4)=(4,-2)$.
  We were quite disappointed to discover that with this definition, \emph{all}
  theorems in the present section would cease to hold.
\end{rmk}

%%%%%%%%%%%%%%%%%%%%%%%%%%%%%%%%%%%%%%%%%%%%%%%%%%%%%%%%%%%%%%%%%%%%%%%%%%%%%%%
\section{Crossings and nestings in set partitions of type $B$}
%%%%%%%%%%%%%%%%%%%%%%%%%%%%%%%%%%%%%%%%%%%%%%%%%%%%%%%%%%%%%%%%%%%%%%%%%%%%%%%
\label{sec:cross-nest-B}

The definition of non-crossing set partitions of type $B$ coincides with the
definition in type $C$, only the combinatorial model for non-nesting set
partitions changes slightly: we define the \Dfn{nesting order} for type $B$ as
$$1<2<\dots<n<0<-n<\dots<-2<-1.$$

\begin{figure}
  \centering
  \setlength{\unitlength}{25pt}
  \begin{picture}(10,2)
    \put( 0,0){\hbox{$1$}}
    \put( 1,0){\hbox{$2$}}
    \put( 2,0){\hbox{$3$}}
    \put( 3,0){\hbox{$4$}}
    \put( 4,0){\hbox{$5$}}
    \put( 5.3,0){\hbox{$0$}}
    \put( 6,0){\hbox{$-5$}}
    \put( 7,0){\hbox{$-4$}}
    \put( 8,0){\hbox{$-3$}}
    \put( 9,0){\hbox{$-2$}}
    \put(10,0){\hbox{$-1$}}
    \qbezier(0.2,0.5)(0.65,1.5)(1.1,0.5)
    \qbezier(1.2,0.5)(3.75,3.5)(6.3,0.5)
    \qbezier(2.2,0.5)(2.65,1.5)(3.1,0.5)	
    \qbezier(3.2,0.5)(4.25,2.5)(5.3,0.5)
    \qbezier(4.2,0.5)(6.8,3.5)(9.4,0.5)	
    \qbezier(5.5,0.5)(6.45,2.5) (7.4,0.5)	
    \qbezier(7.5,0.5)(7.95,1.5) (8.4,0.5)
    \qbezier(9.5,0.5)(9.95,1.5)(10.4,0.5)
  \end{picture}\\
  (a) \\
  \setlength{\unitlength}{25pt}
  \begin{picture}(10,3)
    \put(0,0){\hbox{$1$}}
    \put(1,0){\hbox{$2$}}
    \put(2,0){\hbox{$3$}}
    \put(3,0){\hbox{$4$}}
    \put(4,0){\hbox{$5$}}
    \put(5,0){\hbox{$-1$}}
    \put(6,0){\hbox{$-2$}}
    \put(7,0){\hbox{$-3$}}
    \put(8,0){\hbox{$-4$}}
    \put(9,0){\hbox{$-5$}}
    \qbezier(0.2,0.5)(0.65,1.5)(1.1,0.5)
    \qbezier(1.2,0.5)(2.15,2)(3.1,0.5)
    \qbezier(2.2,0.5)(4.85,3)(7.3,0.5)	
    \qbezier(3.2,0.5)(6.25,4.5)(9.3,0.5)
    \qbezier(4.2,0.5)(6.25,3.5)(8.3,0.5)
    \qbezier(5.4,0.5)(5.85,1.5)(6.3,0.5)
    \qbezier(6.4,0.5)(7.35,2)(8.3,0.5)
  \end{picture}\\
  (b)
  \caption{The nesting (a) and the crossing (b) diagram of a set partition
    of type $B_5$.}
  \label{NNB5}
\end{figure}
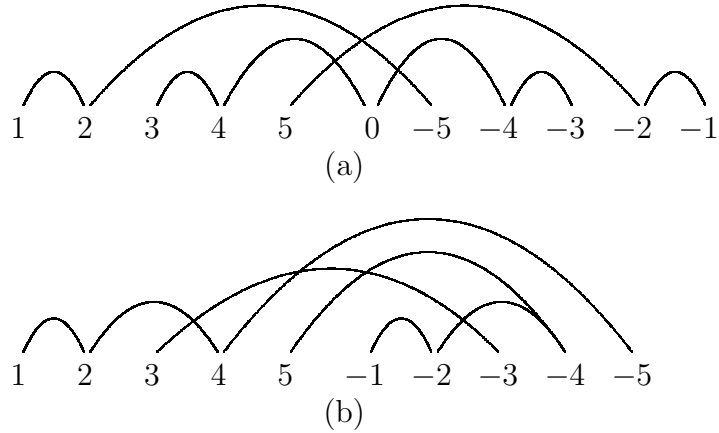

The \Dfn{nesting diagram} of a set partition $\B$ is obtained by placing the
numbers in $[\pm n]\cup 0$ in \emph{nesting order} on a line and then joining
consecutive elements of each block of $\B$ by an arc, where the zero block is
augmented by the number $0$, if present.  See Figure~\ref{NNB5}(a) for an
example.  The number $0$ is neither an opener nor a closer.

These changes are actually dictated by the general, type independent
definitions for non-crossing and non-nesting set partitions.  However, it turns
out that we moreover need to ignore certain crossings and nestings that appear
in the diagrams:

A \Dfn{crossing} is a pair of arcs that crosses in the crossing diagram, except
if both arcs have positive opener and negative closer, and at least one of them
has a closer that is smaller in absolute value than the corresponding opener.
Similarly, a \Dfn{nesting} is a pair of arcs that nests in the nesting diagram,
except if both arcs have positive opener (or begin at $0$) and negative closer
(or end at $0$), and at least one of them has a closer that is smaller in
absolute value than the corresponding opener.

\begin{eg}
  The set partition in Figure~\ref{NNB5}(b) has three crossings: $(3,-3)$
  crosses $(2,4)$, $(4,-5)$, and $(-4,-2)$.  It does not cross $(5,-4)$ by
  definition.

  The set partition in Figure~\ref{NNB5}(a) has three nestings: $(2,-5)$ nests
  $(3,4)$ and $(4,0)$, and $(5,-2)$ nests $(-4,-3)$.  However, $(5,-2)$ does
  not nest $(0,-4)$ by definition.
\end{eg}

With this definition, we have a theorem that is only slightly weaker than in
type $C$:
\begin{thm}
  There is an explicit bijection on set partitions of type $B$, preserving the
  set of openers and the set of closers, and mapping the number of nestings to
  the number of crossings.
\end{thm}

Again, we obtain an analog of Corollary~\ref{cor:unique-A}:
\begin{cor}
  For any opener-closer configuration $(\Set{O},\Set{C}) \subseteq [n] \times
  [n]$, there exists a unique non-crossing set partition $\B$ and a unique
  non-nesting set partition $\B'$, both of type $B_n$, such that
  $$
  \op(\B)=\op(\B')=\Set{O} \quad\text{and}\quad \cl(\B)=\cl(\B')=\Set{C}.
  $$
\end{cor}
\begin{proof}  
  The first two steps of the bijection described in the proof of
  \ref{thm:crossings-C} can be adopted unmodified for the present situation.
  However, it is no longer the case that the notions of nesting and crossing
  coincide for arcs with positive or zero opener and negative or zero closer.

  We remark that there is still exactly one non-nesting way to connect the
  remaining active openers $\{o_1, o_2,\dots,o_k\}$ with their negative
  counterparts, and the number $0$ if $k$ is odd, such that the zero block
  contains $0$ and the symmetry property \eqref{eq:setpartitionsymmetry} is
  satisfied.  For example, the situation for $k=3$ is as follows:
  \setlength{\unitlength}{18pt}
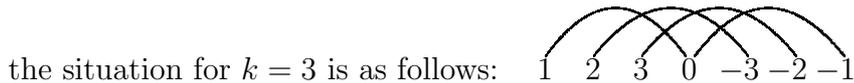
  \begin{picture}(7,2)
    \put(0,0){\hbox{ $1$}}
    \put(1,0){\hbox{ $2$}}
    \put(2,0){\hbox{ $3$}}
    \put(3,0){\hbox{ $0$}}
    \put(4,0){\hbox{$-3$}}
    \put(5,0){\hbox{$-2$}}
    \put(6,0){\hbox{$-1$}}
    \qbezier(0.4,0.5)(1.85,2.5)(3.3,0.5)
    \qbezier(1.4,0.5)(3,2.5)  (4.6,0.5)
    \qbezier(2.4,0.5)(4,2.5)  (5.6,0.5)	
    \qbezier(3.5,0.5)(5,2.5)  (6.6,0.5)
  \end{picture}

  It remains to describe more generally a bijection that maps a type $B$ set
  partition $\B$ with opener-closer configuration $(\Set O, \Set C)=([k],
  \emptyset)$ with $l$ nestings to a type $B$ set partition with $l$ crossings,
  and the same opener-closer configuration.  In fact, we will really map $\B$
  to a type $C$ set partition, such that there are exactly $l$ nestings
  occurring in the set of arcs $(o,c)$ with $o<\abs{c}$.  This is sufficient,
  since for type $C$ set partitions, two arcs $(i,j)$ and $(i',j')$ where both
  $i$ and $i'$ are positive and both $j$ and $j'$ are negative cross if and
  only if they nest.

  If $\B$ does not contain a zero block, the image under the bijection is $\B$
  itself.  Otherwise, suppose that $\B$ consists of arcs 
  $$(o_1, c_1=0), (o_2,c_2), \dots, (o_m, c_m),$$
  with $o_i \leq \abs{c_i}$ for $i>1$, together with their negatives.  We
  assume furthermore that $\abs{c_2} > \abs{c_3} > \dots > \abs{c_m}$, i.e.,
  the closers appear in nesting order.

  Now let $j$ be minimal such that $o_j>\abs{c_{j+1}}$, or, if no such $j$
  exists, set $j:=m$.  We then set
  \begin{equation*}
    (\tilde o_i, \tilde c_i):=
    \begin{cases}
      (o_i, c_{i+1})&\text{for $i<j$}\\
      (o_i, -o_i)&\text{for $i=j$}\\
      (o_i, c_i)&\text{for $i>j$.}
    \end{cases}
  \end{equation*}

  We need to show that the number of nestings among
  $$(\tilde o_1, \tilde c_1), (\tilde o_2,\tilde c_2), \dots, (\tilde o_m,
  \tilde c_m)$$ is the same as in the original set of arcs.  It is sufficient
  to show $\tilde c_{j-1} < \tilde c_j < \tilde c_{j+1}$, i.e., $c_j < -o_j <
  c_{j+1}$, since all other order relations remain unchanged.  The relation
  $o_j < -c_j$ was required for all arcs, and $o_j > -c_{j+1}$ follows from the
  definition of $j$.
\end{proof}

Together with Theorem~\ref{thm:k-crossing-C}, the bijection employed in the
previous proof also shows the following theorem:
\begin{thm}\label{thm:k-crossing-B}
  There is an explicit bijection on set partitions of type $B$, preserving the
  set of openers and the set of closers, and interchanging the cardinalities of
  the maximal crossing and the maximal nesting.
\end{thm}

%%%%%%%%%%%%%%%%%%%%%%%%%%%%%%%%%%%%%%%%%%%%%%%%%%%%%%%%%%%%%%%%%%%%%%%%%%%%%%%
\section{Non-crossing and non-nesting set partitions in type $D$}
%%%%%%%%%%%%%%%%%%%%%%%%%%%%%%%%%%%%%%%%%%%%%%%%%%%%%%%%%%%%%%%%%%%%%%%%%%%%%%%
\label{sec:cross-nest-D}

In type $D$ we do not have any good notion of crossing or nesting, we can only
speak properly about \emph{non-crossing} and \emph{non-nesting} set partitions.

% The \Dfn{nesting order} and the \Dfn{nesting diagram} are defined as in type
% $B$ and $C$.  However, the meaning of \Dfn{nesting} is not as intuitive as in
% the other types.  In type $D_n$,

\begin{figure}
  \setlength{\unitlength}{25pt}
  \begin{picture}(9,2.5)
    \put(0,0){\hbox{$1$}}
    \put(1,0){\hbox{$2$}}
    \put(2,0){\hbox{$3$}}
    \put(3,0){\hbox{$4$}}
    \put(4,0){\hbox{$5$}}
    \put(5,0){\hbox{$-1$}}
    \put(6,0){\hbox{$-2$}}
    \put(7,0){\hbox{$-3$}}
    \put(8,0){\hbox{$-4$}}
    \put(9,0){\hbox{$-5$}}
    \qbezier(0.2,0.5)(4.25,4.5)(8.3,0.5)
    \qbezier(1.2,0.5)(2.65,2.5)(4.1,0.5)
    \qbezier(2.2,0.5)(7.25,3.5)(9.3,0.5)
    \qbezier(3.2,0.5)(4.25,1.5)(5.3,0.5)
    \qbezier(4.2,0.5)(5.75,2.5)(7.3,0.5)
    \qbezier(6.4,0.5)(7.85,2.5)(9.3,0.5)
  \end{picture}
  \setlength{\unitlength}{25pt}
  \begin{picture}(9,2.5)
    \put(0,0){\hbox{$1$}}
    \put(1,0){\hbox{$2$}}
    \put(2,0){\hbox{$3$}}
    \put(3,0){\hbox{$4$}}
    \put(4,0){\hbox{$5$}}
    \put(5,0){\hbox{$-1$}}
    \put(6,0){\hbox{$-2$}}
    \put(7,0){\hbox{$-3$}}
    \put(8,0){\hbox{$-4$}}
    \put(9,0){\hbox{$-5$}}
    \qbezier(0.2,0.5)(4.25,3.5)(8.3,0.5)
    \qbezier(1.2,0.5)(6.75,3.5)(9.3,0.5)
    \qbezier(2.2,0.5)(3.15,1.5)(4.1,0.5)
    \qbezier(3.2,0.5)(4.25,1.5)(5.3,0.5)
    \qbezier(4.2,0.5)(5.25,1.5)(6.3,0.5)
    \qbezier(7.4,0.5)(8.35,1.5)(9.3,0.5)
  \end{picture}
  \caption{Two non-crossing set partition of type $D_5$. Both are obtained from
    each other by interchanging $5$ and $-5$.}
  \label{NCD5}
\end{figure}
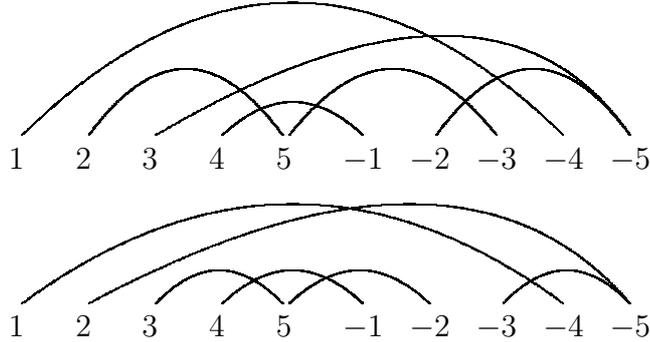

A combinatorial model for non-crossing set partition of type $D_n$ was given by
Christos A. Athanasiadis and Vic Reiner in~\cite{AthanasiadisReiner2004}.
% the \Dfn{crossing diagram} is obtained by placing the numbers from $1$ to
% $n-1$ and $-1$ to $-(n-1)$ in this order on the vertices of a $2(n-1)$-gon,
% and $n$ and $-n$ in its middle.
For our purposes it is easier to use a different description of the same model:
let $\B$ be a set partition of type $D_n$ and let $\{(i_1,-j_1), \dots,
(i_k,-j_k)\}$ for positive $i_l,j_l < n$ be the ordered set of arcs in $\B$
starting in $\{1,\ldots,n-1\}$ and ending in its negative. $\B$ is called
\Dfn{non-crossing} if
\begin{itemize}
\item[(i)] $(i,-i)$ is an arc in $\B$ implies $i = n$,
\end{itemize}
and if it is non-crossing in the sense of type $C_n$ with the following
exceptions:
\begin{itemize}
\item[(ii)] arcs in $\B$ containing $n$ \emph{must} cross all arcs $(i_l,-j_l)$
  for $l > k/2$,
\item[(iii)] arcs in $\B$ containing $-n$ \emph{must} cross all arcs $(i_l,-j_l)$
  for $l \leq k/2$,
\item[(iv)] two arcs in $\B$ containing $n$ and $-n$ \emph{may} cross.
\end{itemize}
Here, (i) is equivalent to say that if $\B$ contains a zero block $B_0$ then $n
\in B_0$ and observe that (i) together with the non-crossing property of
$\{(i_1,-j_1),\ldots,(i_k,-j_k)\}$ imply that $k/2 \in \mathbb{N}$, see
Figure~\ref{NCD5} for an example.

\begin{rmk}
  All conditions hold for a set partition $\B$ if and only if they hold for the
  set partition obtained from $\B$ by interchanging $n$ and $-n$.
\end{rmk}
                
\begin{figure}
  \setlength{\unitlength}{25pt}
  \begin{picture}(9,2)
    \put(0,0){\hbox{$1$}}
    \put(1,0){\hbox{$2$}}
    \put(2,0){\hbox{$3$}}
    \put(3,0){\hbox{$4$}}
    \put(4,0){\hbox{$5$}}
    \put(5,0){\hbox{$-5$}}
    \put(6,0){\hbox{$-4$}}
    \put(7,0){\hbox{$-3$}}
    \put(8,0){\hbox{$-2$}}
    \put(9,0){\hbox{$-1$}}
    \qbezier(0.2,0.5)(2.75,2.5)(5.3,0.5)
    \qbezier(1.2,0.5)(2.75,1.5)(4.1,0.5)
    \qbezier(2.2,0.5)(4.25,2.5)(6.3,0.5)
    \qbezier(3.2,0.5)(5.25,2.5)(7.3,0.5)
    \qbezier(4.2,0.5)(6.75,2.5)(9.3,0.5)
    \qbezier(5.4,0.5)(6.85,1.5)(8.3,0.5)
  \end{picture}
  \setlength{\unitlength}{25pt}
  \begin{picture}(9,2)
    \put(0,0){\hbox{$1$}}
    \put(1,0){\hbox{$2$}}
    \put(2,0){\hbox{$3$}}
    \put(3,0){\hbox{$4$}}
    \put(4,0){\hbox{$5$}}
    \put(5,0){\hbox{$-5$}}
    \put(6,0){\hbox{$-4$}}
    \put(7,0){\hbox{$-3$}}
    \put(8,0){\hbox{$-2$}}
    \put(9,0){\hbox{$-1$}}
    \qbezier(0.2,0.5)(2.15,2.5)(4.1,0.5)
    \qbezier(1.2,0.5)(3.25,2.5)(5.3,0.5)
    \qbezier(2.2,0.5)(4.25,2.5)(6.3,0.5)
    \qbezier(3.2,0.5)(5.25,2.5)(7.3,0.5)
    \qbezier(4.2,0.5)(6.25,2.5)(8.3,0.5)
    \qbezier(5.2,0.5)(7.25,2.5)(9.3,0.5)
  \end{picture}
  \caption{Two non-nesting set partition of type $D_5$. Both are obtained from
    each other by interchanging $5$ and $-5$.}
  \label{NND5}
\end{figure}
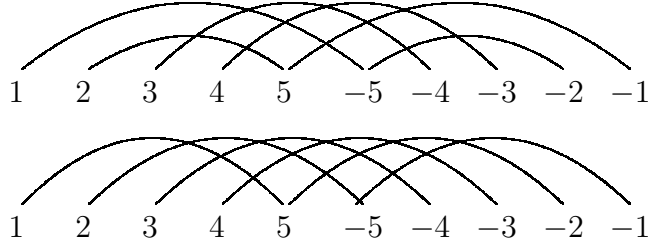

A set partition of type $D_n$ is called \Dfn{non-nesting} if it is non-nesting
in the sense of \cite{Athanasiadis1998}. This translates to our notation as
follows: let $\B$ be a set partition of type $D_n$. Then $\B$ is called
non-nesting if
\begin{itemize}
\item[(i)] $(i,-i)$ is an arc in $\B$ implies $i = n$,
\end{itemize}
and if it is non-nesting in the sense of type $C_n$ with the following
exceptions:
\begin{itemize}
\item[(ii)] arcs $(i,-n)$ and $(j,n)$ for positive $i<j<n$ in $\B$ are allowed
  to nest, as do
\item[(iii)] arcs $(i,-j)$ and $(n,-n)$ for positive $k < i,j < n$ in $\B$
  where $(k,n)$ is another arc in $\B$ (which exists by the definition of set
  partitions in type $D_n$).
\end{itemize}
Again, (i) means that if $B_0 \in \B$ is a zero block then $n \in B_0$.  (ii)
and (iii) come from the fact that the positive roots $e_i + e_n$ and $e_j -
e_n$ for $i \leq j$ are comparable in the root poset of type $C_n$ but are not
comparable in the root poset of type $D_n$, see Figure~\ref{NND5} for an
example.  As for non-crossing set partitions in type $D_n$, all conditions hold
if and only if they hold for the set partition obtained by interchanging $n$
and $-n$.
 
\begin{prop}\label{pr:NCD} 
  Let $(\Set{O},\Set{C}) \subseteq [n]$ be an opener-closer configuration. Then
  there exists a non-crossing set partition $\B$ of type $D_n$ with $\op(\B) =
  \Set{O}$ and $\cl(\B) = \Set{C}$ if and only if
  \begin{eqnarray} 
    |\Set{O}| - |\Set{C}| \text{ is even or } n \in \Set{O},\Set{C}. \label{eq:NCD}
  \end{eqnarray} 
  Moreover, there exist exactly two non-crossing set partitions of type $D_n$
  having this opener-closer configuration if both conditions hold, otherwise,
  it is unique.
\end{prop}
                        
\begin{proof} 
  Suppose that $|\Set{O}| - |\Set{C}|$ is odd. Then the conditions to be
  non-crossing imply that we must have a zero block and therefore, $n$ must be
  an opener.  On the other hand, the definition of set partitions of type $D_n$
  implies that $n$ must be a closer. Thereby, condition (\ref{eq:NCD}) is
  necessary.  For the proof of the proposition we distinguish three cases:
                                
  Case 1: $|\Set{O}| = |\Set{C}|$. Then by the definition of opener-closer
  configurations, $n \notin O$ and the unique construction is the same as in
  the first step of the proof of Theorem~\ref{thm:crossings-C}.
                                
  Case 2: $|\Set{O}| - |\Set{C}|$ is odd. Then by (\ref{eq:NCD}), $n$ is both
  opener and closer. For $\Set{C} \setminus \{n\}$ the construction is the same
  as in Case 1. Now, there is an odd number of positive openers smaller than
  $n$ left. Connect the closer in $n$ to the unique opener in the middle as
  well as the opener in $-n$ to its negative. Connect $n$ and $-n$.  Finally
  connect the remaining openers with there negative counterparts as closers
  such that they are non-crossing.
                                
  Case 3: $|\Set{O}| - |\Set{C}| > 0$ is even. For $\Set{C} \setminus \{n\}$
  the construction is again as in type $C_n$. Now, there is an even number of
  positive openers left. If $n$ is a closer but not an opener, then there is an
  odd number of positive openers smaller than $n$ left. Connect the closer in
  $n$ to the unique opener in the middle as well as the opener in $-n$ to its
  negative. If $n$ is a closer and also an opener then there is an even number
  of positive openers smaller than $n$ left. Connect the closer in $n$ to one
  of the two openers in the middle and the opener in $n$ to the negative of the
  other and also connect $-n$ to their negatives. This gives the two
  possibilities in this case and observe that both are obtained from each other
  by interchanging $n$ and $-n$. Finally connect the remaining openers with
  there negative counterparts as closers such that they are non-crossing.
\end{proof}

As in types $A, B$ and $C$, the analogue proposition holds also for non-nesting
set partitions of type $D_n$:

\begin{prop}\label{pr:NND} 
  Let $(\Set{O},\Set{C}) \subseteq [n]$ be an opener-closer configuration. Then
  there exists a non-nesting set partition $\B$ of type $D_n$ with $\op(\B) =
  \Set{O}$ and $\cl(\B) = \Set{C}$ if and only if
  \begin{eqnarray} 
    |\Set{O}| - |\Set{C}| \text{ is even or } n \in \Set{O},\Set{C}. \label{eq:NND}
  \end{eqnarray} 
  Furthermore, there exist exactly two non-nesting set partitions of type $D_n$
  having this opener-closer configuration if both conditions hold, otherwise,
  it is unique.
\end{prop}

\begin{proof} 
  The proof that condition \eqref{eq:NND} is necessary is analogous to the
  proof in the non-crossing case.

  Recall that a set partition of type $D_n$ is non-nesting if it is non-nesting
  in the sense of type $C_n$ except for arcs of the forms
  \begin{itemize}
  \item[(i)] arcs $(i,-n)$ and $(j,n)$ for positive $i<j<n$,
  \item[(ii)] arcs $(i,-j)$ and $(n,-n)$ for positive $k < i,j < n$ where
    $(k,n)$ is another arc (which exists if $(n,-n)$ is an arc),
  \end{itemize} 
  and observe that in both cases, $n$ is both an opener and a
  closer. Therefore, the construction is exactly the same as in type $C_n$
  otherwise. We now prove the remaining two cases:
                                
  Case 1: $|\Set{O}| - |\Set{C}|$ is odd.  The unique possibility is to connect
  $n$ and $-n$ and all others in the same way as in type $C_n$. All nesting
  arcs in this case are of the form (ii).
                                
  Case 2: $|\Set{O}| - |\Set{C}|$ is even. In this case, we have two
  possibilities: the first is to connect closers and openers as in type $C_n$
  without creating any nestings. The second is to connect the closers in
  $\Set{C} \setminus \{n\}$ as above to the associated openers, then we connect
  $-n$ to the first active opener and $n$ to the associated negative
  closer. The remaining positive openers and their associated negative closers
  are finally connected such that they are non-nesting. Observe All nesting
  arcs in this case are of the form (i). Observe also that possibilities 1 and
  2 are obtained from each other by interchanging $n$ and $-n$.
\end{proof}

%%%%%%%%%%%%%%%%%%%%%%%%%%%%%%%%%%%%%%%%%%%%%%%%%%%%%%%%%%%%%%%%%%%%%%%%%%%%%%%
\section{$k$-crossing and $k$-nesting set partitions of type $C$}
%%%%%%%%%%%%%%%%%%%%%%%%%%%%%%%%%%%%%%%%%%%%%%%%%%%%%%%%%%%%%%%%%%%%%%%%%%%%%%%
\label{sec:k-crossing-C}

In this section we prove Theorem~\ref{thm:k-crossing-C}, which states that the
cardinalities of the maximal crossing and the maximal nesting of type $C$ set
partitions are equidistributed.

The rough idea of our bijection is as follows: we first show how to render a
$C_n$ set partition in the language of $0$-$1$-fillings of a certain polyomino,
as depicted in Figure~\ref{fig:growth}(a).  We will do this in such a way that
maximal nestings correspond to north-east chains of ones of maximal length.  

Interpreting this filling as a growth diagram in the sense of Sergey Fomin and
Tom Roby~\cite{Fomin1986,Roby1991,Fomin1994,Fomin1995} enables us to define a
transformation on the filling that maps -- technicalities aside -- the length
of the longest north-east chain to the length of the longest south-east chain.
This filling can then again be interpreted as a $C_n$ set partition, where
south-east chains of maximal length correspond to maximal crossings.  Many
variants of the transformation involved are described in Christian
Krattenthaler's article~\cite{Krattenthaler2006}, we will employ yet another
(slight) variation.

Let us now give a detailed description of the objects involved: the
\Dfn{nesting polyomino} for type $C_n$ set partitions is the polyomino
consisting of $n$ columns of height $2n-1, 2n-2, \dots, n$, arranged in this
order.  We label the columns $1,2,\dots, n$ and the rows from top to bottom
$2,3,\dots,n,-n,\dots, \dots, -2, -1$, as in Figure~\ref{fig:growth}(a).  Thus,
every box of the polyomino corresponds to an arc with positive opener, that may
be present in a nesting diagram: an arc $(i, j)$ corresponds to the cell in
column $i$, row $j$.

We encode a type $C_n$ set partition by placing ones into those boxes that
correspond to arcs, and zeroes into the other boxes, as in
Figure~\ref{fig:growth}(a).  (For convenience, zeros are not shown and ones are
indicated by crosses.  We ignore the integer partitions labelling the top-right
corners for the moment.)  A $0$-$1$-filling of the nesting polyomino
corresponds to a set-partition if and only if
\begin{enumerate}
\item there is at most one non-zero entry in each row and each column
\item the restriction of the filling to the rows $-1,-2,\dots, -n$ is symmetric
  with respect to the diagonal as indicated in the figure, and
\item there is at most one non-zero entry on this diagonal.
\end{enumerate}

The \Dfn{crossing polyomino} for type $C_n$ set partitions is a polyomino of
the same shape as the nesting polyomino.  We label the columns $1,2,\dots n$ as
before.  However, we now label the columns the rows from top to bottom
$2,3,\dots,n,-1,-2, \dots, -n$, as in Figure~\ref{fig:growth}(b).  We find that
a $0$-$1$-filling of the crossing polyomino corresponds to a set-partition
under the same conditions as before, with the difference that the symmetry axis
now runs south-east instead of north-east.

\def\dr{\POS[];[d]**\dir{-},[r]**\dir{-}}
\def\r{\POS[];[r]**\dir{-}}
\def\d{\POS[];[d]**\dir{-}}
\def\dd{\POS[];[d]**\dir{.}}
\def\ddr{\POS[];[d]**\dir{.},[r]**\dir{-}}
\def\e{\emptyset}
\def\x{\mbox{\Large$\times$}}
\def\y{\POS[];[dr]**{}?*{\x}}
\def\put#1{\POS[];[dr]**{}?*{\mbox{#1}}}
\def\di{\POS[];[-5,5]**\dir{.}}
\def\id{\POS[];[5,5]**\dir{.}}
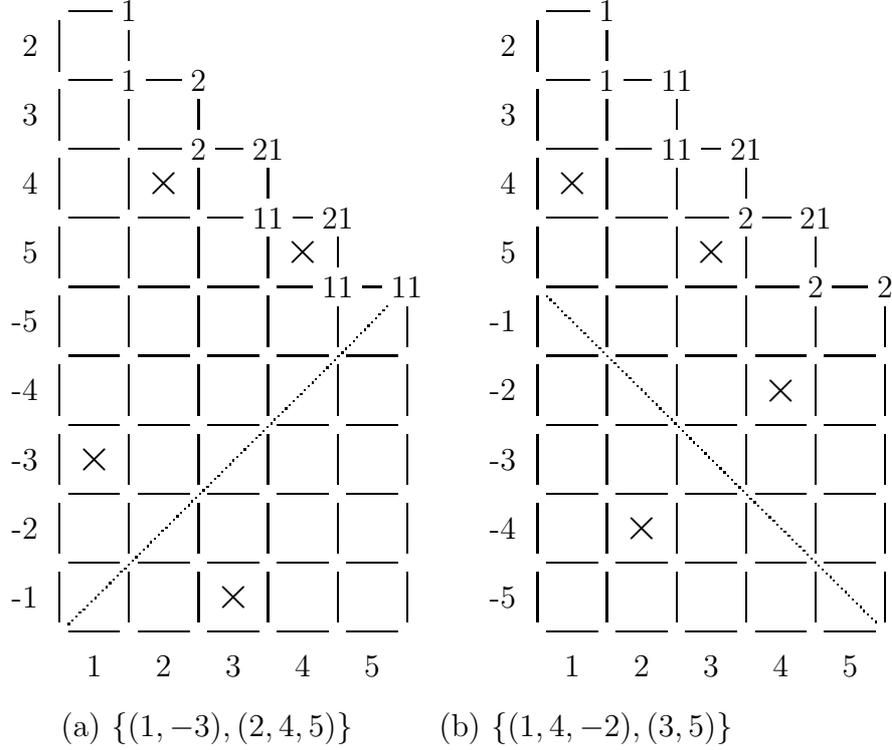
\begin{figure}[h]
  \begin{tabular}{cc}
    \begin{xy}<0pt,126pt>;<10pt,126pt>:
      \xymatrix@!=2pt{%
        {}\put{ 2}&{}\dr  &1 \d   &       &       &       &    \\
        {}\put{ 3}&{}\dr  &1 \dr  &2 \d   &       &       &    \\
        {}\put{ 4}&{}\dr  &{}\dr\y&2 \dr  &21\d   &       &    \\
        {}\put{ 5}&{}\dr  &{}\dr  &{}\dr  &11\dr\y&21\d   &    \\
        {}\put{-5}&{}\dr  &{}\dr  &{}\dr  &{}\dr  &11\dr  &11\d\\
        {}\put{-4}&{}\dr  &{}\dr  &{}\dr  &{}\dr  &{}\dr  &{}\d\\
        {}\put{-3}&{}\dr\y&{}\dr  &{}\dr  &{}\dr  &{}\dr  &{}\d\\
        {}\put{-2}&{}\dr  &{}\dr  &{}\dr  &{}\dr  &{}\dr  &{}\d\\
        {}\put{-1}&{}\dr  &{}\dr  &{}\dr\y&{}\dr  &{}\dr  &{}\d\\
        &{}\r\put{1}\di&{}\r\put{2}&{}\r\put{3}&{}\r\put{4}&{}\r\put{5}&\\ 
        &{}     &{}     &{}     &{}     &{}      &{} }  
    \end{xy}
    &
    \begin{xy}<0pt,126pt>;<10pt,126pt>:
      \xymatrix@!=2pt{%
        {}\put{ 2}&{}\dr   &1 \d   &       &       &       &    \\
        {}\put{ 3}&{}\dr   &1 \dr  &11\d   &       &       &    \\
        {}\put{ 4}&{}\dr\y &{}\dr  &11\dr  &21\d   &       &    \\
        {}\put{ 5}&{}\dr   &{}\dr  &{}\dr\y&2 \dr  &21\d   &    \\
        {}\put{-1}&{}\dr\id&{}\dr  &{}\dr  &{}\dr  &2 \dr  &2\d\\
        {}\put{-2}&{}\dr   &{}\dr  &{}\dr  &{}\dr\y&{}\dr  &{}\d\\
        {}\put{-3}&{}\dr   &{}\dr  &{}\dr  &{}\dr  &{}\dr  &{}\d\\
        {}\put{-4}&{}\dr   &{}\dr\y&{}\dr  &{}\dr  &{}\dr  &{}\d\\
        {}\put{-5}&{}\dr   &{}\dr  &{}\dr  &{}\dr  &{}\dr  &{}\d\\
        &{}\r\put{1}&{}\r\put{2}&{}\r\put{3}&{}\r\put{4}&{}\r\put{5}&\\ 
        &{}     &{}     &{}     &{}     &{}      &{} }  
    \end{xy}\\
    (a) $\{(1,-3),(2,4,5)\}$ & \hspace{-2.5cm}(b) $\{(1,4,-2),(3,5)\}$
  \end{tabular}
  \caption{growth diagrams for type $C_5$ set partitions}
  \label{fig:growth}
\end{figure}

A \Dfn{north-east chain} of length $k$ is a sequence of $k$ non-zero entries in
a filling of a nesting polyomino, such that every entry is strictly to the
right and strictly above the preceding entry in the sequence.  Similarly, a
\Dfn{south-east chain} of length $k$ is a sequence of $k$ non-zero entries in a
filling of a crossing polyomino, such that every entry is strictly to the right
and strictly below the preceding entry in the sequence.  Furthermore, we
require that the smallest rectangle containing all entries of the sequence is
completely contained in the polyomino.  We remark that this condition is
trivially satisfied for north-east chains in fillings of nesting polyominoes.

\begin{lem}\label{lem:preserve-length}
  A longest north-east chain in a $0$-$1$-filling of the nesting polyomino
  corresponds to a maximal nesting in the corresponding set partition.
  Similarly, a longest south-east chain in a $0$-$1$-filling of the crossing
  polyomino corresponds to a maximal crossing in the corresponding
  set-partition.
\end{lem}
\begin{proof}
  The statement for the nesting polyomino is trivial.  For the crossing
  polyomino we have to show that if a maximal crossing involves arcs with
  negative opener and negative closer, there is another maximal crossing that
  involves arcs with positive openers only.

  We note that a maximal crossing cannot contain an arc with positive opener
  and positive closer, and an arc with negative opener and negative closer
  simultaneously.

  Thus, by symmetry, if a maximal crossing $(o_1,c_1), (o_2,c_2), \dots, (o_k,
  c_k)$ involves an arc with negative opener and negative closer, the set of
  arcs $(-c_1,-o_1), (-c_2,-o_2), \dots, (-c_k, -o_k)$ is also a maximal
  crossing, with all openers positive.
\end{proof}

We can now explain the significance of the integer partitions labelling the
top-right corners along the border of the upper half of the polyominoes in
Figure~\ref{fig:growth}: the sum of the first $k$ parts of each of these
partitions is just the maximal cardinality of a union of $k$ \emph{north-east}
chains in the rectangular region of the polyomino to the left and below the
corner.  Moreover, the sum of the first $k$ parts of the conjugate partition
equals the maximal cardinality of a union of $k$ \emph{south-east} chains.

In particular the first part of every partition equals the length of the
longest north-east chain in the region under consideration, while the first
part of its conjugate equals the length of the longest south-east chain.  As an
aside, we remark that the sum of the parts gives the number of non-zero entries
within this region.

The following proposition is a consequence of the general theory of growth
diagrams:
\begin{prop}\label{prop:growth-bijection}
  A $0$-$1$-filling of a \emph{nesting polyomino} corresponds bijectively to a
  type $C_n$ set partition, if and only if
  \begin{enumerate}
  \item\label{item:vacillating} the sequence of integer partitions $(\lambda_1,
    \lambda_1,\dots,\lambda_{2n-1})$ labelling the top-right corners along the
    border of the upper half of the polyomino, when read from top to bottom, is
    \Dfn{vacillating}.  I.e., for all $k$ we have
    \begin{itemize}
    \item $\lambda_{2k-1}=\lambda_{2k}$, or $\lambda_{2k-1}$ is obtained from
      $\lambda_{2k}$ by adding one to some part, and
    \item $\lambda_{2k+1}=\lambda_{2k}$, or $\lambda_{2k+1}$ is obtained from
      $\lambda_{2k}$ by adding one to some part.
    \end{itemize}
  \item\label{item:diagonal} the bottom most integer partition contains at most
    one column of odd length, i.e., the conjugate partition has at most one odd
    part.
  \end{enumerate}

  A $0$-$1$-filling of a \emph{crossing polyomino} corresponds bijectively to a
  type $C_n$ set partition, if and only if
  \begin{enumerate}
  \item[(1')] the sequence of integer partitions $(\lambda_1,
    \lambda_1,\dots,\lambda_{2n-1})$ labelling the top-right corners along the
    border of the upper half of the polyomino, when read from top to bottom, is
    vacillating,
  \item[(2')]\label{item:diagonal'} and the bottom most integer partition
    contains at most one one odd part.
  \end{enumerate}
\end{prop}
\begin{proof}
  Let us first consider nesting polyominoes.  Suppose we are given a filling
  that corresponds to a type $C_n$ set partition.  In the preceding paragraphs,
  we already described how to obtain the integer partititions labelling the
  top-right corners along the border of the upper half of the polyomino in
  Figure~\ref{fig:growth}.  The vacillating condition~\eqref{item:vacillating}
  is satisfied, since there is at most one non-zero entry in every column and
  every row.  Since the filling restricted to rows $-n$ to $-1$ is symmetric,
  and there is at most one non-zero entry on the diagonal,
  Condition~\eqref{item:diagonal} is satisfied by, for example \cite[Exercise
  7.28]{EC2}.

  We now have to show how to recover a filling given only the sequence of
  partitions.  Using the \lq local backward rules\rq\ $(B1)$--$(B6)$ as
  defined, for example in \cite[Section~2]{Krattenthaler2006}, we can recover
  the entries in rows $2$ to $n$ of the $0$-$1$-filling, as well as a sequence
  of integer partitions labelling the top-right corners of row $-n$.  It
  remains to find out how to label the top-right corners of column $n$, so we
  can also determine the filling in rows $-n$ to $-1$.

  It is well known (eg., \cite[Corollary~7.13.6]{EC2} that the filling of a
  square growth diagram is symmetric with respect to its main diagonal, if and
  only if the sequence of partitions labelling the top-right corners along the
  top-most row is the same as the sequence of partitions labelling the the
  top-right corners along the right-most column.  Thus, we have no choice but
  to label the top-right corners of column $n$ with the sequence of partitions
  we just computed for the top-right corners of row $-n$.

  The proof for crossing polyominoes is very similar, we only have to explain
  how to obtain the filling of rows $-1$ to $-n$, given the sequence of integer
  partitions labelling the top-right corners of row $-1$.  Let $Q$ be the
  (partial) standard Young tableau corresponding to this sequence of
  partitions.
  
  By \cite[Corollary~A1.2.11]{EC2} we know that reflecting the filling
  restricted to rows $-1$ to $-n$ about a vertical axis corresponds to the
  following transformation:
  \begin{itemize}
  \item the integer partitions labelling the top-right corners of column $n$
    are all conjugated, whereas
  \item the integer partitions labelling the top-right corners of row $-1$ are
    obtained by evacuating the (partial) standard Young tableau, and then
    transposing the corresponding partitions.
  \end{itemize}
  In particular, since evacuation is an involution, the filling restricted to
  rows $-1$ to $-n$ is symmetric with respect to the diagonal indicated in
  Figure~\ref{fig:growth}(b) if and only if the sequence of partitions labelling
  the top-right corners of column $n$ correspond to the (partial) standard
  Young tableau obtained by evacuating $Q$.
\end{proof}

It is now obvious how to construct the desired bijection demonstrating
\ref{thm:k-crossing-C}:
\begin{proof}[Proof of Theorem~\ref{thm:k-crossing-C}]
  \begin{enumerate}
  \item given a $0$-$1$-filling of a nesting polyomino, compute the sequence of
    integer partitions labelling the top-right corners along the border of its
    upper half,
  \item label the top-right corners along the border of the upper half of a
    crossing polyomino with the conjugate partitions
  \item using \ref{prop:growth-bijection} compute the $0$-$1$-filling
    corresponding to the labelling.
  \end{enumerate}
\end{proof}

We have to remark that the bijection presented above is not an involution.
Furthermore, it does not exchange the crossing and the nesting numbers.  As a
small example, consider the $C_4$ partition $\{1,4\},\{2,-3\}$, which is
non-nesting, has four crossings, and the cardinality of the maximal crossing is
two.  Its crossing polyomino is mapped to the nesting polyomino of the $C_4$
partition $\{1, -3\},\{2,4\}$, which has two nestings, two crossings.  Of
course, by construction of the bijection, the cardinality of the maximal
nesting is two, also.
%%%%%%%%%%%%%%%%%%%%%%%%%%%%%%%%%%%%%%%%%%%%%%%%%%%%%%%%%%%%%%%%%%%%%%%%%%%%%%%
\section{$k$-triangulations in set partitions of type $C$ and $k$-fans of
  symmetric Dyck paths}
%%%%%%%%%%%%%%%%%%%%%%%%%%%%%%%%%%%%%%%%%%%%%%%%%%%%%%%%%%%%%%%%%%%%%%%%%%%%%%%
\label{sec:k-triangulations-C}

In this section we want to deduce a conjecture due to Daniel Soll and Volkmar
Welker~\cite{SollWelker2006}, using the same methods as in the previous
section.  Namely, we consider consider generalized triangulations of the
$2n$-gon that are invariant under rotation of $180\degree$, and such that at
most $k$ diagonals are allowed to cross mutually.  Daniel Soll and Volkmar
Welker then conjectured that the number of such triangulations that are
maximal, i.e., where one cannot add any diagonal without introducing a $k+1$
crossing, coincides with the number of \Dfn{fans} of $k$ Dyck paths that are
symmetric with respect to a vertical axis.  We start with the precise
definitions:

Consider the $2n$-gon with vertices labelled clockwise from
$$1,2,\dots,n,-1,-2,\dots,-n.$$
Let $\omega$ be a set of diagonals that is invariant under rotation of
$180\degree$.  I.e. if the diagonal $\{i,j\}$ is present, then the diagonal
$\{-i,-j\}$ must be present, too.  Obviously, every set partition of type $C_n$
(or of type $B_n$) can be regarded as such a subset of diagonals, by including
exactly those diagonals that connect labels in the crossing diagram of the set
partition.

A subset of $k$ diagonals of $\omega$ that mutually cross in the relative
interior of the polygon is a \Dfn{$k$-crossing}.  We remark that two diagonals
$\{i,j\}$ and $\{k,l\}$ \Dfn{cross} exactly if $i<k<j<l$ in the \Dfn{crossing
  order} $1<2<\dots<n<-1<-2<\dots<-n$.  Thus, the notion of crossing we
describe here agrees with the notion of crossing in
Section~\ref{sec:cross-nest-C}.  To avoid a misconception that distracted the
authors for some time, we stress the fact that $\{i,j\}$ and $\{i',j'\}$ need
not cross even if $\{i,j\}$ and $\{-i',-j'\}$ do.

We now encode $\omega$ as in the previous section by a $0$-$1$-filling of the
crossing polyomino, placing ones into those boxes that correspond to diagonals.
Note that the number of non-zero entries above and including the main diagonal
in the filling is just the number of diagonals in $\omega$.  Again, we have
that a longest south-east chain in the filling of the crossing polyomino
corresponds to a maximal crossing in the corresponding set-partition.  If
$\omega$ is maximal, that is, adding any diagonal increases the cardinality of
the maximal crossing, and its maximal crossing has cardinality $k$, we call
$\omega$ a \Dfn{type $C_n$ $k$-triangulation}.

The second kind of objects under consideration are \Dfn{symmetric fans of $k$
  non-intersecting Dyck paths}.  For our purposes it is best to define them as
families of paths in the nesting polyomino: the paths start in the boxes
labelled $(1,2), (2,3),\dots,(k,k+1)$, take unit south or west steps, and end
on the main diagonal.  Furthermore, we insist that they are non-intersecting.

Let us call a $0$-$1$-filling of a nesting polyomino \Dfn{maximal}, when
replacing a zero by a one in any box increases the length of the longest
north-east chain.  It is easy to construct a bijection between such fans and
maximal fillings of the nesting polyomino whose longest north east chain has
length $k$: we simply put a one in every box that a path enters, and zeroes
elsewhere.

We can now state the main theorem of this section:
\begin{thm}[Conjecture~13 of \cite{SollWelker2006}]\label{thm:k-triangulation}
  The number of type $C_n$ $k$-triangulations coincides with the number of
  symmetric fans of $k$ non-intersecting Dyck paths.  Equivalently, the number
  of maximal $0$-$1$ fillings of a nesting polyomino whose length of the
  longest north-east chain equals $k$ coincides with the number of maximal
  $0$-$1$ fillings of a crossing polyomino whose length of the longest
  south-east chain equals $k$.
\end{thm}

The corresponding theorem for type $A$ was discovered and proved by Jakob
Jonsson~\cite{Jonsson2005}.  A (nearly) bijective proof very similar to ours
was given by Christian Krattenthaler in~\cite{Krattenthaler2006}.  A simple
bijection for the case of $2$-triangulations was recently given by Sergi
Elizalde in~\cite{Elizalde2006}.

The main difference to the previous section is that there will now be several
non-zero entries in most of the rows and columns of the polyominoes.  Thus, we
have to use a variant of the bijection, for \Dfn{arbitrary fillings} of
polyominoes with non-negative integers, constructed in the previous section,
and deduce Theorem~\ref{thm:k-triangulation} inductively thereafter.

%% define vertical strip

\begin{prop}\label{prop:growth-bijection-semi}
  An arbitrary filling of a \emph{nesting polyomino} corresponds bijectively to
  a sequence of integer partitions $(\lambda_1,
  \lambda_1,\dots,\lambda_{2n-1})$ labelling the top-right corners along the
  border of the upper half of the polyomino, when read from top to bottom, if
  and only if for all $k$ we have
  \begin{itemize}
  \item $\lambda_{2k-1}=\lambda_{2k}$, or $\lambda_{2k-1}$ is obtained from
    $\lambda_{2k}$ by adding a horizontal strip, and
  \item $\lambda_{2k+1}=\lambda_{2k}$, or $\lambda_{2k+1}$ is obtained from
    $\lambda_{2k}$ by adding a horizontal strip.
  \end{itemize}

  An arbitrary filling of a \emph{crossing polyomino} corresponds bijectively
  to a sequence of integer partitions $(\lambda_1,
  \lambda_1,\dots,\lambda_{2n-1})$ labelling the top-right corners along the
  border of the upper half of the polyomino, when read from top to bottom, if
  and only if for all $k$ we have
  \begin{itemize}
  \item $\lambda_{2k-1}=\lambda_{2k}$, or $\lambda_{2k-1}$ is obtained from
    $\lambda_{2k}$ by adding a vertical strip, and
  \item $\lambda_{2k+1}=\lambda_{2k}$, or $\lambda_{2k+1}$ is obtained from
    $\lambda_{2k}$ by adding a vertical strip.
  \end{itemize}
\end{prop}
\begin{proof}
  The general procedure is as in the proof of
  Proposition~\ref{prop:growth-bijection}.  However, we now have to use
  different \lq local backward\rq\ rules, since we are dealing with arbitrary
  fillings.  Namely, in the case of nesting polyominoes, we use the rules $(B^1
  0)$--$(B^1 2)$ of \cite[Section~4.1]{Krattenthaler2006}, whereas in the case
  of crossing polyominoes we use the rules $(B^4 0)$--$(B^4 2)$ of
  \cite[Section~4.4]{Krattenthaler2006}.
\end{proof}

\begin{proof}[Proof of Theorem~\ref{thm:k-triangulation}]
  \begin{enumerate}
  \item given an arbitrary filling of a nesting polyomino, compute the sequence
    of integer partitions labelling the top-right corners along the border of
    its upper half,
  \item label the top-right corners along the border of the upper half of a
    crossing polyomino with the conjugate partitions
  \item using \ref{prop:growth-bijection-semi} compute the filling
    corresponding to the labelling.
  \end{enumerate}
  Note that this filling will in general not be a $0$-$1$-filling.  However, we
  remark that the sum of all the entries in the filling of the nesting
  polyomino and in the filling of the crossing polyomino will be the same.

  Still, we can deduce the statement of the theorem.  Let $N(m,l)$ be the set
  of fillings of nesting polyominoes with length of longest north-east chain
  equal to $k$, sum of all entries equal to $m$ and $l$ non-zero entries.
  Similarly, let $C(m,l)$ be the set of fillings of crossing polyominoes with
  length of longest north-east chain equal to $k$, sum of all entries equal to
  $m$ and $l$ non-zero entries.

  We will prove that $\abs{N(m,l)}=\abs{C(m,l)}$, the particular case
  $\abs{N(m,m)}=\abs{C(m,m)}$ is exactly the statement of the theorem.  When
  $m$ equals one, the cardinalities of the fillings coincide by
  Proposition~\ref{thm:k-triangulation} of the previous section, since in this
  case we actually have at most one non-zero entry in every row and every
  column.  The statement for arbitrary $m$ and $l=1$ follows trivially.

  More generally, if $\abs{N(l,l)}=\abs{C(l,l)}$, it follows that
  $\abs{N(m,l)}=\abs{C(m,l)}$ for $m>l$: let $m=m_1+m_2+\dots+m_l$ be a
  composition of $m$ into $l$ parts.  Then every possibility of replacing the
  $l$ ones in a filling in $N(l,l)$ by $m_1,m_2,\dots,m_l$ corresponds
  bijectively to a possibility of replacing the $l$ ones in a filling in
  $C(l,l)$ by $m_1,m_2,\dots,m_l$.

  We know already that the number of arbitrary fillings of the nesting
  polyomino and the number of arbitrary fillings of the crossing polyomino with
  sum of all entries equal to $m$ coincide.  This number equals
  $$\abs{N(m,1)}+\abs{N(m,2)}+\dots+\abs{N(m,m)},$$
  but also
  $$\abs{C(m,1)}+\abs{C(m,2)}+\dots+\abs{C(m,m)}.$$
  By induction, we know that $\abs{N(m,l)}=\abs{C(m,l)}$ for $l<m$.  Therefore,
  $\abs{N(m,l)}$ and $\abs{C(m,l)}$ must coincide, too.
\end{proof}

\bibliographystyle{hamsplain} %amsplain, amsalpha
%\bibliography{math.bib}
\providecommand{\cocoa} {\mbox{\rm C\kern-.13em o\kern-.07em C\kern-.13em
  o\kern-.15em A}}
\providecommand{\bysame}{\leavevmode\hbox to3em{\hrulefill}\thinspace}
\providecommand{\href}[2]{#2}

\end{document}